\documentclass[11pt]{amsart}

\usepackage{amsmath,amssymb,amscd,amsfonts,verbatim}
\usepackage{paralist}
\usepackage[mathscr]{eucal}
\usepackage{hyperref}
\usepackage{color}
\usepackage{pdfsync}

\newtheorem{thm}{Theorem}[section]
\newtheorem{lem}[thm]{Lemma}

\newtheorem{prop}[thm]{Proposition}
\newtheorem{cor}[thm]{Corollary}

\newtheorem{assu-nota}[thm]{Assumption--Notation}
\theoremstyle{remark}

\newcommand{\pp}{\mathbb P}

\newcommand{\OO}{\mathcal O}

\DeclareMathOperator{\Pic}{Pic}

\DeclareMathOperator{\Alb}{Alb}

\numberwithin{equation}{section}

\title{A characterization of the symmetric square of a curve}
\author{Margarida Mendes Lopes, Rita Pardini and Gian Pietro Pirola}

\begin{document}
\begin{abstract} \medskip
In this paper  a new intrinsic geometric characterization of the symmetric square of a curve and  of the ordinary  product of two curves is given. More precisely it is shown that the existence on a  surface of general type $S$ of irregularity $q$ of an effective divisor $D$ having self-intersection $D^2>0$ and arithmetic genus  $q$ implies that $S$ is either birational to a product of curves or to the second symmetric  product of a curve. 


\end{abstract}
\maketitle
%\tableofcontents
\section{Introduction}

In this paper we give a new intrinsic geometric characterization of the symmetric square of a curve and  of the ordinary  product of two curves. \bigskip

 The symmetric products $S^k(C)$ of a curve $C$ of genus $g>0$ give,  together with the ordinary products,
the simplest examples of irregular varieties.  To give some perspective, we recall
that the cohomology ring of $S^k(C)$ approximates the cohomology ring of the Jacobian $J(C)$ of $C.$ 
In particular one has (\cite{Mac}) $H^1(C,\mathbb Z)\equiv H^1(J(C),\mathbb Z)\equiv H^1(S^k(C),\mathbb Z)$ and, for $k>1$, 
$$H^2(S^k(C),\mathbb Z) \equiv H^2(J(C),\mathbb Z)\oplus \mathbb Z.$$ 

Notice  that  the Torelli-type theorem for $k<g$ (see \cite{ran2}) provides an equivalence between curves and symmetric products.
\bigskip

In view of their  simple cohomological structure, one expects the symmetric products
 to  have a very important place in the classification of irregular varieties. Some progress has been done using Fourier-Mukai transform and generic vanishing in the case $k=g-1$ by C. Hacon (\cite{Hacon}), giving  a cohomological characterization of the theta divisor of  a principally polarized abelian variety (PPAV). 
  \medskip
  
 For $1<k<g-2$ there is  a famous conjecture  due to O. Debarre (\cite{deb2}) that claims that the symmetric products $S^k(C)$ and the Fano surface $F$ of the lines of a smooth cubic $3$-fold are the only varieties giving the minimal cohomological class  of a PPAV  $(A,\theta)$.
The Debarre conjecture is proved for  $q=4,$ by Z. Ran (\cite{ran}), when $A=J(X)$ is a Jacobian of a curve by Debarre (\cite{deb2}),  and when  $A$ is the  intermediate Jacobian of a generic cubic $3$-fold by A. H\"oring (\cite{horing}).

%\noindent {\bf Conjecture 1} {\cite ddebarre 95} \\ {\em Let   $(A,\theta)$ be a PPAV of dimension $q.$ Let  $X$ be a subvariety of $A$ of dimension $1<k<q-1.$ Assume that the cohomological class of $X$ is minimal, that is
% $\theta^{q-k}/(q-k)!.$  \begin{enumerate} \item either  $A$ is the jacobian of a curve $C$ of genus $q$$X=W_k$ or $X=-W_k$ where $W_k$ is the Abel-Jacobi image of the simmetric product $W_k=AJ(C_k)$  of the curve $C.$ \item  or $ q=5$ and $k=2,$ $A=J(V)$ is the intermediate Jacobina of a smooth cubic threefold $V$ and $X$ is isomorphic to the Abel Jacobi image, $X=AJ(F)$ or $X=-AJ(F),$ of the Fano surface of the lines of $V$ .\bigskip\end{enumerate}}
\bigskip 

Coming back to surfaces, it is a basic problem to give some geometrical or cohomological  characterization of the  second symmetric product of a curve.
%Let $\Omega^1_X$ be the cotangent bundle of a surface $X$  of irregularity
 %$q=\dim H^{0}(\Omega^1_X)>2.$
 
  In the spirit of surface classification theory,
a conjecture, not equivalent, and in some  sense even stronger than  Debarre's, is the following:
\bigskip

\noindent {\bf Conjecture: }\\
\noindent {\em The only minimal surfaces $X$ of irregularity $q>2$ with $H^{0}(\Omega^2_X) \equiv \bigwedge^2 H^{0}(\Omega^1_X),$  are the  symmetric products $S^2(C)$ and the Fano surfaces $F$ of the lines of  a smooth cubic $3$-fold.}

\bigskip

%We remark that the Fano surfaces  of the lines $F$ of a smooth cubic $3-$fold  are the only known examples where $H^2(X,\mathbb Z)\equiv \wedge^2 H^1(X,\mathbb Z).$ \\

The  above conjecture was proven for  $q=3$ in \cite{HaconPardini} and independently in \cite{Pirola}. 
The proof  of \cite{HaconPardini}  uses  in a crucial way the  fact  that the image of the Albanese map is a divisor
 and
it seems  difficult to generalize. 
%The Fano surface example also shows that some care is necessary.
The proof of  \cite{Pirola}
is based on the geometric characterization of $S^2(C)$ given
in \cite{ccm}  using  the geometry of families of curves.
Namely,  in  \cite{ccm}   $S^2(C)$ is  proven to be the only minimal
algebraic surface  with irregularity $q$ covered by curves of genus $q$ and self-intersection $1$.

\bigskip
Analysing the curves of small genus on a surfaces of general type, we discovered a surprisingly precise geometric characterization of the  second symmetric product (and of the ordinary product).

   \begin{thm}\label{main} Let $S$ be a smooth surface of general type with irregularity $q$ containing a $1$-connected divisor  $D$ such that $p_a(D)=q$ and $D^2>0$. Then the minimal model of $S$ is either:
   \begin{itemize} \item[(a)] the product of two curves of genus $g_1, g_2\geq 2$ (and $g_1+g_2=q$) or 
   \item[(b)] the symmetric product  $S^2(C)$ where $C$ is a smooth curve of genus $q$ (and $C^2=1$).
   
   \end{itemize}
   Furthermore, if $D$ is $2$-connected, only the second case occurs.
\end{thm}

This result is in some sense very atypical of the theory of algebraic surfaces, because we obtain a complete classification of the surface from the existence of a single divisor with certain properties.    The only similar instance   we know of  is the characterization of rational surfaces from the existence of a smooth rational curve with positive self-intersection.
\bigskip 

The proof of Theorem \ref{main} consists of two main steps, 
that we describe in the symmetric product case.
First we  reduce to the case when the curve $C$ is smooth.  This step requires  a very careful numerical analysis  of the effective divisors  contained in  $C$. Then we observe that the Albanese variety $\Alb(X)$  of the surface  is isomorphic to the Jacobian $J(C)$ of $C$ and,  combining the Brill-Noether theory on $C$ with the generic vanishing theorem of Green and Lazarsfeld,   
we show that the curve $C$ moves in a positive dimensional family. This allows us to use the  results in \cite{ccm}  and complete the proof. 
%We remark that the second step combines the generic vanishing and Brill-Noether theory of the curve $C$. 
\medskip

It appears that a sort of weak Brill-Noether theory can be performed for line bundles on irregular surfaces. We will come back to these topics in a forthcoming paper (cf. \cite{future}).

\bigskip

\noindent{\bf Notation:}  %$\sim$ denotes numerical equivalence of divisors. %%For a surface $S$ $q$ will denote... 
A  {\em surface} is a smooth projective complex surface. The {\em irregularity} of a surface $S$, often denoted by $q$ or $q(S)$, is $h^1(\OO_S)=h^0(\Omega^1_S)$ and the {\em geometric genus} $p_g(S)$ is $h^0(K_S)=h^2(\OO_S)$.
  A {\em curve} on a surface is a nonzero effective divisor. A {\em fibration} of genus $b$ of  a surface $S$ is a map $f\colon S\to B$ with connected fibers, where $B$ is a smooth curve of genus $b$; $f$ is  {\em relatively minimal} if its fibers  contain  no $-1$-curve.    A {\em quasi-bundle} is a fibration such that all the smooth fibers are isomorphic and  the  singular fibers of are multiples of smooth curves.

 \section{Auxiliary facts }
In this section we collect some auxiliary facts.
We start by  recalling the well known equality for reducible curves on  smooth surfaces:
\medskip
  
$\bullet$ If a curve $D$ decomposes as $D=A+B$ where $A,B>0$  and $AB=m$ then $p_a(A)+p_A(B)+m-1=p_a(D)$. 
\smallskip
 
 The following results are needed in the sequel: 
  \begin{lem}\label{trans}    Let $S$ be a  surface and  let $f\colon S\to B$ a relatively minimal fibration of genus $\geq 2$ with general fibre of genus $\geq 2$.  If $C$ is  a section  of $f$ (i.e. an irreducible curve  $C$ such that $CF=1$), then:
  \begin{enumerate}
  \item  $C^2\leq 0$;
  \item if  $C^2=0$,  then  $f$ is a quasi-bundle and there exists $m>0$ such that $mC$ is a fiber  of a quasi-bundle fibration of $S$.
  \end{enumerate}
    \end{lem}
   \begin{proof}     (i) Let us notice first that  any section $C$  of $f$ is smooth of genus $b$. Denote by $F$ the general fiber of $f$; by   Arakelov's theorem (cf. \cite[Theorem 3.1]{Serrano}) the relative canonical class  $K_{S/B}\sim K_S-(2b-2)F$  is nef,  hence $K_SC\ge (2b-2)$ and, by the adjunction formula, $C^2\le0$. 

(ii): follows by  Theorem 3.2 in \cite{Serrano}.
  \end{proof}

    \begin{lem}\label{subcurve} Let $D$ be a 1-connected curve on a surface. Then every subcurve $A<D$ satisfies $p_a(A)\leq p_a(D)$.
  \end{lem}
   \begin{proof} Set $B:=D-A$ and let  $m:=AB$. Then $p_a(A)+p_a(B)+m-1=p_a(D)$.    Suppose for  contradiction that $p_a(A)>p_a(D)$.  Then we obtain $p_a(B)+m-1<0$, i.e. $m<1-p_a(B)$. Since $1-p_a(B)=h^0(B,\OO_B)-h^1(B, \OO_B)$ we conclude that $m<h^0(B,\OO_B)$.  This contradicts  \cite[ Lemma 1.4]{KML}, that states that any subcurve $D'$ of a $1$-connected curve $D$ satisfying $D'(D-D')=b$ verifies $h^0(D',\OO_{D'})\leq b$ .
  
   \end{proof} 
  \begin{lem} \label{pa}
  Let $S$ be a surface of general type and let $D$ be an irreducible curve of $S$. If $D^2>0$, then:
  \begin{enumerate}
  \item $p_a(D)\ge 2$;
  \item if $p_a(D)=2$, then the minimal model $T$ of $S$ has $K^2_T=1$, $q(T)=0$.
  \end{enumerate}
  \end{lem}
  \begin{proof} Let $\eta\colon S\to T$ be the morphism onto the minimal model, so that $K_S:=\eta^*K_T+E$, where $E$ is the exceptional divisor. Since $D^2>0$, the curve $D$ is not contracted by $\eta$, hence $K_SD\ge \eta^*K_TD>0$ and the adjunction formula gives immediately $p_a(D)\ge 2$, with equality holding if and only if $D^2=\eta^*K_T D=1$. 
  Hence, if $p_a(D)=2$ the index theorem gives $K^2_T=1$ and,  since by \cite{deb1} irregular surfaces have $K^2\ge 2p_g$,  $T$ is regular.  \end{proof}
  
      \section{Curves with $p_a=q$   }\label{sec:pq}

      %Now we start the proof of Theorem \ref{main}.  
     As a preparation for the proof of Theorem \ref{main}, in this section we make a detailed numerical analysis of the following situation:
     \begin{itemize}
     \item $S$ is a surface of general type with irregularity $q$;
     \item $D$ is a $1$-connected curve  of $S$ such that $p_a(D)=q$ and $D^2>0$.
     \end{itemize}
     \medskip
     
     %Observe that the previous assumptions imply that $q>0$: indeed it is well known that 
  \begin{prop}\label{product} Let $S$ be a surface of general type   with irregularity $q$ and let  $D$ % Kodaira dimension \geq 0
 be a  $1$-connected curve of $S$     such that $p_a(D)=q$ and $D^2>0$.
If $D=A+B$, where $A,B$ are curves such that  $AB=1$ and $p_a(A)\geq 1$, $p_a(B)\geq 1$,  then $S$ is  birational to the product of two curves.
 \end{prop}

\begin{proof} 
  The equality $AB=1$ implies that both $A$ and $B$ are 1-connected   (see \cite[Lemma A.4]{cfm}) and that  $q=p_a(D)=p_a(A)+p_a(B)$.  Since, by assumption,  $p_a(A)\geq 1$, $p_a(B)\geq 1$, one has  also $p_a(A)<q$ and $p_a(B)<q$. So  both $A^2\leq 0$, $B^2\leq 0$ (cf. \cite{ra}, \cite[Remark 6.8]{cat} and also \cite{q6}).  By the  hypothesis $D^2=A^2+2+B^2>0$, at least one of the inequalities is in fact  an equality. Suppose then that $B^2=0$. Then (cf. ibidem) there exists a fibration $f\colon S\to E$ where $E$ is a curve of geometric genus $g(E)\geq q-p_a(B)$ and such that $mB$ is a fibre $F$ of $f$ for some $m>0$.

Since $AB=1$ and $B$ is nef,  there is a unique irreducible curve $\theta\leq A$ such that $\theta B\neq 0$.  Of course $\theta$  is not contained in a fibre of $f$ and,  by Lemma \ref{subcurve}, $p_a(\theta)\leq p_a(A)$. From $p_a(A)=q-p_a(B)$ we obtain $p_a(\theta)\leq g(E)$.   Since $f$ induces a surjective morphism $\theta\to E$ of degree $m$ and $p_a(\theta)\leq g(E)$, we conclude, by the Hurwitz formula,  that  $g(\theta)=p_a(\theta)=g(E)$ and, in addition, $m=1$ or  $g(E)=1$ and $\theta \to E$ is unramified. In the latter case,  if $m> 1$  we have a contradiction because the  existence of a multiple fibre of $f$  means that  the cover is ramified.
So $m=1$ and  the fibration $f\colon S\to E$ satisfies $g(E)+g(F)=q$. By \cite[Lemme]{appendix},  $S$ is birational to the product  $F\times E$.
\end{proof} 
\begin{comment}
\medskip

The above proposition, that corresponds to the first case in  Theorem \ref{main}, has the following immediate:

\smallskip

  \begin{cor}\label{corproduct} Let $S$ be a surface of general type of irregularity $q$ 
  %Kod\geq 0
   and   assume that $S$  is not   birational to the product of two curves. Let $D$ be a 1-connected curve of $S$ such that $p_a(D)=q$ and $D^2>0$.
  Then every decomposition of $D=A+B$ where  $A,B>0$ and $AB=1$ satisfies $p_a(A)=0$ or $p_a(B)=0$.   
 \end{cor}
\end{comment}

Next we show that, when $S$   is not   birational to the product of two curves, we can assume that $D$ is 2-connected.

  \begin{lem}\label{1con} Let $\tilde S$ be an irregular surface of general type
  %Kod \geq 0
   that   is not   birational to the product of two curves. Let $D$ be a 1-connected curve of $\tilde S$ such that $p_a(D)=q$, $D^2>0$. Then there is a birational morphism $p\colon\tilde S\to S_0$, where $S_0$ is a smooth surface such that $p(D)$ contains a 2-connected curve $D_0$ satisfying $D_0^2>0$ and $p_a(D_0)=q$. 
%  Furthermore if $D$ is nef  then there is a birational morphism $p\colon\tilde S\to S$, where $S$ is a smooth surface such that $f(D)$ is a 2-connected curve  satisfying $f(D)^2=D^2>0$ and $p_a(f(D))=q$.
 
 \end{lem}
\begin{proof} 
  If $D$ is 2-connected there is nothing to prove. Otherwise we have a decomposition $D=A+B$ where $AB=1$ and, by  Lemma \ref{product},  $p_a(A)=0$, $p_a(B)=q$.  By  \cite[Lemma A.4]{cfm}, $A$ and $B$ are $1$-connected.  Since $p_a(A)=0$,  $A$ is contracted by the Albanese map and  so $A^2<0$. 
  
  If $A^2<-1$, then from  $D^2>0$ and $AB=1$ we must have $B^2>0$ and we consider now  the curve $B$.    If $B$ is 2-connected we take $D_0:=B$.  If $B$ is not  2-connected, then  it has  a decomposition $B=A_1+A_2$ with $A_1A_2=1$. As above we can suppose $p_a(A_1)=0$ and  $p_a(A_2)=q$ and we can now restart the reasoning.  
  
   If  $A^2=-1$,  then  $A$  contains an irreducible $(-1)$-curve $E$  such that $ED=0$.  If $p\colon S\to S_1$ is the blow down of $E$ and $D_1:=p(D)$, then we have $D=p^*D_1$, hence $D_1$ is $1$-connected, $D_1^2=D^2>0$ and $p_a(D_1)=p_a(D)=q$. Hence we may replace $D$ and $S$ by $D_1$ and $S_1$ and, if $D_1$ is not $2$-connected, repeat the previous step.
      
   Since $D$ has a finite number of components, the process described above  must stop and so in the end we get a surface $S_0$ birational to $S$ contained a  $2$-connected $D_0$ as wished.
%Suppose now that $D$ is nef. Then  if $D$ is not 2-connected  a decomposition $D=A+B$ where $AB=1$ and $p_a(A)=0$, $p_a(B)=q$ must be such that $A^2=-1$ and therefore contracting $A$ we obtain $p(D)$ with the same properties and again nef.
\end{proof}

\section{The $2$-connected case}

Here  we examine  the situation of \S \ref{sec:pq} under the additional assumption that   $D$ be $2$-connected.

  \begin{lem}\label{irr} Let $S$ be an irregular surface of general type and  let $D$ be a 2-connected curve of $S$ such that $p_a(D)=q$ and $D^2>0$.
 %Kod\geq 0 and p_g\geq 0
   Then:
 \begin{enumerate}
 \item $D$ is contained in the fixed part of $|K_S+D|$;
 \item $D$ is smooth irreducible;
 \item $h^0(S, D)=1$;
 \item the Albanese image of $S$ is a surface.
%\item $h^0(D, K_S)<q$. 
 \end{enumerate}
 \end{lem}
  \begin{proof} (i) Since $D^2>0$ and, because $D$ is 1-connected,  $h^0(D,\OO_D)=1$,  one has $h^1(S, K_S+D)=0$. Hence the cokernel of the restriction map  $H^0(S,K_S+D)\to H^0(D, \omega_D)$ is $H^1(S, K_S)$ and therefore by $p_a(D)=q$,  the image of $r$ must be zero. 
  \smallskip
  
 (ii) Assume by contradiction that $D$ is not smooth. Then $D$ has multiple points.  Since $D$  is contained in the fixed part of $|K_S+D|$, any multiple point  $P$ of $D$ is a base point of $|K_S+D|$  and so, by  \cite[Thm. 3.1]{adjoint}, $D$ is not 2-connected. This contradicts the hypothesis. Finally, $D$, being smooth and $1$-connected, is necessarily irreducible.
  \smallskip  
  
  (iii) Assume that $h^0(S,D)>1$.  Then the irreducible curve $D$ moves. Since, by hypothesis, $p_g(S)>0$, we conclude that  $D$     is not a fixed component of $|K_S+D|$, contradicting (i). 
  \smallskip
  
  (iv)  Since $D$ is irreducible we have $q\ge 3$ by Lemma \ref{pa}.  Assume for  contradiction that  the Albanese image of $S$ is a curve $E$.  Then $E$ is a smooth curve of genus $q$ and we have thus a fibration $f\colon S\to E$. Since the smooth curve $D$ satisfies  $D^2>0$, $D$ is not contained in any fibre $F$ of $f$, and so $f$ induces a degree $m$  morphism $f|_D\colon D\to E$, where $m:=DF$. 
 Since $g(D)=g(E)=q$ and $q\geq 3$, we have $m=1$ by the Hurwitz formula.
 
Now notice that we can assume that $S$ is minimal, and thus $f$ is relatively minimal. In fact every $(-1)$-curve $\theta$ is contained in the fibres of $S$ and thus, since  $DF=1$, the image of $D$ in the minimal model of $S$ is still a smooth curve with geometric genus $q$.     Then $D$ is a section of $f$, but this contradicts Lemma \ref{trans}.

  %Since $S$ is of general type $g(F)\geq 2$ and so by Arakelov's theorem (cf. \cite{appendix}) $\omega_{S/B} D>0$.  From $g(D)=q$ and $D^2>0$ one has $K_SD<2q-2$. Then, since  $\omega_{S/B}\sim K_S-(2q-2)F$, we have a contradiction. 
So the Albanese image of $S$ must be  a surface.
 % Finally  since $d:=D^2>0$, $K_SD=2q-2-d<2q-2$ and so by the Riemann- Roch theorem and  Clifford's lemma  $h^0(D,K_S)< q$. 
  \end{proof}

  \begin{cor}\label{num}  Let $S$ be an irregular surface of general type that is not birational to a product of curves  and  let $D$ be a nef 2-connected curve of $S$ such that $p_a(D)=q$ and $D^2>0$.  Assume also that  there is no $(-1)$-curve $\theta$ such that $D\theta=0$. Then any curve $C$ numerically equivalent to $D$ is smooth irreducible.

 \end{cor} 
 \begin{proof}   Since $D$ is nef,  also $C$ is nef.  Then it is well known that $C$ nef and big implies that  $C$ is  $1$-connected (see \cite[Lemma 2.6]{adjoint}).  
If $C$ is not 2-connected there is  a decomposition $C=A+B$ where $AB=1$ and, by Proposition  \ref{product},  $p_a(A)=0$, $p_a(B)=q$.  As in Lemma \ref{1con} one has $A^2<0$ and since $C$ is nef one must have $A^2=-1$, but this contradicts  the hypothesis on $D$.   Hence $C$ is 2-connected and so,  by Lemma   \ref{irr},  $C$ is smooth irreducible.
   
   \end{proof}   

 %\section{The proof of Theorem \ref{main}}
 \begin{prop}\label{Cmoves} Let $S$ be a surface of general type  with irregularity $q$ and let $D$ be  a smooth irreducible curve $D$ such that $d:=D^2>0$ and $g(D)=q$. 
 Then the set $\{L\in \Pic^0(S)|h^0(D+\eta)>0\}$ has dimension $\ge \min\{q,d\}$.
 
 \end{prop}
\begin{proof} By Lemma \ref{irr},   the Albanese image of $S$ is a surface.
 So, by the generic vanishing theorem of Green-Lazarsfeld (\cite{GL1},\cite{GL2}), the set $V^1(S)=\{\eta\in \Pic^0(S)|h^1(\eta)>0\}$ is the union of finitely many translates of proper abelian subvarieties of $A$. 
Since $D^2>0$, the map $\Pic^0(S)\to \Pic^0(D)$ is injective and thus an isomorphism, hence we can identify $\Pic^0(S)$ and $\Pic^0(D)$.

 Denote by $W_d(D)\subset \Pic^0(D)$ the image of the natural map $S^d(D)\to \Pic^0(D)$ defined by $\Delta\mapsto [\OO_D(\Delta -D)]$.  If $d\ge q$ then $W_d(D)=\Pic^0(D)$, otherwise $W_d(D)$ is  $d$-dimensional. Then  $W_d(D)$   generates $\Pic^0(D)$  and it is irreducible, hence   it cannot be contained in $V^1(S)$. %, because, as stated above,  $V^1(S)$ the union of finitely many translates of proper abelian subvarieties.   
 
 Note that for $\eta\notin V^1(S)$ the restriction sequence $0\to\eta\to\eta+D\to (\eta+D)|_D\to 0$ gives an isomorphism $H^0(\eta+D)\cong H^0((\eta+D)|_D)$. Hence for every $\eta\in W_d(D)\setminus V^1(S)$ (and, by semicontinuity,  for every $\eta\in W_d(D)$) we have $h^0(D+\eta)>0$.
 \end{proof}

 \section{Proof of the main result}
 
 \begin{proof}[Proof of Theorem \ref{main}]
 
  By Proposition \ref{product}, if there exists a decomposition $D=A+B$, with $A,B>0$, $AB=1$ and $p_a(A),p_a(B)>0$, then $S$ is birational to a product of curves, namely we have case (a).
  
  So assume that no such decomposition exists. Then, up to replacing $S$ by a surface in the same birational class, by Proposition \ref{1con} we may assume that $D$ is $2$-connected, hence smooth and irreducible by Lemma \ref{irr}.

Write $d:=D^2$. By Proposition \ref{Cmoves}, there exist a $d$-dimensional system $\mathcal C$ of curves numerically equivalent to $D$.   Since we can obviously assume that  there is no $(-1)$-curve $A$ such that $DA=0$,  any curve $C$ of $\mathcal C$ is  smooth by Lemma  \ref{num}. The Jacobian of  every  curve of $\mathcal C$ is isomorphic to $\Pic^0(D)$, hence the smooth elements of $\mathcal C$ are all isomorphic to $D$.

 If $d>1$, by  \cite[Lemma 2.2.1]{guerra} (cf. also \cite[\S0]{ccm})  $S$ is not of general type, against the assumptions.  So we have $d=1$ and the result follows by \cite[Thm. 0.20]{ccm}.
 \end{proof}   

\bigskip

\noindent {\em Acknowledgments.}  {\small The first author is a member of the Center for Mathematical
Analysis, Geometry and Dynamical Systems (IST/UTL).  The second and the third author are members of G.N.S.A.G.A.--I.N.d.A.M. This research was partially supported by  FCT (Portugal) through program POCTI/FEDER and Project 
PTDC/MAT/099275/2008 and by the italian  PRIN 2008 project  {\it Geometria delle
variet\`a algebriche e dei loro spazi di moduli}.}

\bigskip

\begin{minipage}{13.0cm}
\parbox[t]{6.5cm}{Margarida Mendes Lopes\\
Departamento de  Matem\'atica\\
Instituto Superior T\'ecnico\\
Universidade T{\'e}cnica de Lisboa\\
Av.~Rovisco Pais\\
1049-001 Lisboa, PORTUGAL\\
mmlopes@math.ist.utl.pt
 } \hfill
\parbox[t]{5.5cm}{Rita Pardini\\
Dipartimento di Matematica\\
Universit\`a di Pisa\\
Largo B. Pontecorvo, 5\\
56127 Pisa, Italy\\
pardini@dm.unipi.it}

\vskip1.0truecm

\parbox[t]{5.5cm}{Gian Pietro Pirola\\
Dipartimento di Matematica\\
Universit\`a di Pavia\\
Via Ferrata, 1 \\
 27100 Pavia, Italy\\
\email{gianpietro.pirola@unipv.it}}
\end{minipage}

\end{document}